\documentclass[11pt,reqno]{amsart}
\usepackage{graphicx}
\usepackage{verbatim}
\usepackage{textcomp}
\usepackage{amssymb}
\usepackage{cite}
\usepackage{amsmath}
\usepackage{latexsym}
\usepackage{amscd}
\usepackage{amsthm}
\usepackage{mathrsfs}
\usepackage{xypic}
\usepackage{bm}
\usepackage{url}
\usepackage{hyperref}

\newtheorem{thm}{Theorem}[section]
\newtheorem{corr}[thm]{Corollary}
\newtheorem{lem}[thm]{Lemma}

\theoremstyle{definition}

\theoremstyle{remark}

\numberwithin{equation}{section}
\begin{document}
\title[Gradient estimates of a nonlinear elliptic equation]{Liouville theorems and gradient estimates of a nonlinear elliptic equation for the V-Laplacian}
\author{Yike Jia}
\address{School of Mathematics and Statistics\\
Xinyang Normal University\\
Xinyang, 464000, P. R. China.} \email{1136904169$@$qq.com}
\maketitle
\begin{quotation}
\noindent{\bf Abstract.}~In this paper we establish gradient estimates for positive solutions to the nonlinear elliptic equation $$\Delta_{V}u^{m}+\mu(x)u+p(x)u^{\alpha}=0 , \quad m>1$$on any smooth metric measure space whose $k$-Bakry-\'{E}mery curvature is bounded from below by
$-(k-1)K$ with $K \geq 0$. Additionally, we obtain related Liouville theorems and Harnack inequalities. We partially extend conclusions of Wang, when $V=0$, $\mu=0$ the equation becomes $\Delta u^{m}+p(x)u^{\alpha}=0$. And $V=f$, $\mu=c, p=0 $, the equation becomes $\Delta_{f}u^{m}+cu=0 $.
\end{quotation}

\section{Introduction}
In $1981$, Gidas-Spruck \cite {B81}studied the equation $$\Delta u+hu^{\alpha}=0$$ on a complete manifold with nonnegative Ricci curvature where $1 \leq \alpha \leq \frac{n+2}{n-2}$ and $h(x)\geq0$, if $u(x)$ is a nonnegative solution to the equation, then $u(x)\equiv0$. In $2010$, Yang \cite {Y10} studied the equation
$$\Delta u+cu^{-\alpha}=0$$ on a noncompact complete Riemannian manifold, where $\alpha > 0$ and c are
two realconstants.In 2015, Wang \cite{L16} deduced gradient estimates and Liouville type theorem for
positive solutions to the equation$$\Delta_{f}u^{m}+cu=0$$ on smooth measure space with $m$-Bakry-\'{E}mery curvature bounded by $Ric_{f,m} \geq -(m-1)K$, where $K \geq 0$.Also in 2015, Wang \cite{W17} studied the equation $$\Delta u^{m}+\lambda(x)u^{l}=0$$ on a complete Riemannian manifold without
boundary. Where $m > 1$ and $l$ are real numbers, and $\lambda(x) \in C^{2}(M^{n})$.

An important generalization of Laplacian is the V-Laplacian
$$\Delta_V\star=\Delta\star+\langle V, \nabla\star\rangle$$
on a Riemannian manifold $(M^{n},g)$ where $V$ is a smooth field on $M$. A smooth
function $u$ is called a $V$-harmonic function if it is a solution to$$\Delta_{V}u=0$$
 $V$-harmonic function is a special case of $V$-harmonic maps introduced by Chen,
Jost, and Wang. As in\cite {G09,Q12}, we define the $k$-Bakry-\'{E}mery Ricci tensor
$$Ric_V^k=Ric-\frac{1}{2}L_Vg-\frac{1}{k-n}V^\flat\otimes V^\flat,$$
$$Ric_V=Ric-\frac{1}{2}L_Vg,$$
where Ric is the Ricci tensor of $M$, $L$ is the Lie derivative, $\flat$ maps a vector
field to its dual $1$-form isomorphically, and $k \geq n$ is a constant. Throughout
this paper, we agree that $k = n$ if and only if $V=0$.
From the Bochner formula for $\Delta_{V}$
$$\frac{1}{2}\Delta_V|\nabla u|^2=|\operatorname{Hess} u|^2+\operatorname{Ric}_V(\nabla u,\nabla u)+\langle\nabla u,\nabla\Delta_Vu\rangle.$$

Inspired by the works by\cite {B81,Y10,L16,W17,G20}we investigate the nonlinear elliptic equation
\begin{equation}\label{1Int1}
\Delta_{V}u^{m}+\mu(x)u+p(x)u^{\alpha}=0, \quad m>1,
\end{equation}
on a complete Riemannian manifold with $k$-Bakry-\'{E}mery Ricci curvature bounded below, where
$m > 1$ and $\alpha$ are real numbers, and $p(x),\mu(x) \in C^{2}(M^{n}). $

\begin{thm}\label{Int-1} Let $(M^{n},g)$ be a complete Rimannian mainfold without boundary. Assume that $B_{2R}$ is a geodesic ball of radius $2R$ centered at $q\in M$ and $Ric_V^k(B_{2R}) \geq -K$
with $K \geq 0$. Additionally, suppose there are two positive numbers $\delta_{1},\delta_{2}$ such that
$|p(x,t)|\leq\delta_{1},|\mu(x,t)|\leq\delta_{2}$. Furthermore, assume that $\mu p\geq0$ , $|\nabla p|^{2} \leq |p|^{2}$ and $|\nabla \mu|^{2} \leq |\mu|^{2}$. Let $u(x)$ be a positive solution to the equation \eqref{1Int1}
and $v=\frac{m}{m-1}u^{m-1}$.

$( 1) $ Assume that $\alpha \geq 1$ , $\beta > 0$ then

\begin{equation}
\sup_{B_q(R)}\frac{|\nabla v|^2}{v}\leq\frac{C_4(m-1)}{m+1}[\frac{1}{R^2}(1+\sqrt{K}R)+2K+1]\sup_{x\in M^{n}}v+H_1.
\end{equation}

$( 2) $ Assume that $\alpha < 1$ , $\beta > 0$  then
\begin{equation}\sup_{B_q(R)}\frac{|\nabla v|^2}{v}\leq\frac{C_4(m-1)}{m+1}[\frac{1}{R^2}(1+\sqrt{K}R)+2K+1]\sup_{x\in M^{n}}v+H_2.
\end{equation}
Where $C_{4}$ is a constant depending on $k$,
$$\aligned
 H_{1}&=\frac{(m-1)(k-1)}{m(m+1)}|\frac{2(m+1)}{n-1}+(m-2\alpha+1)|\delta_{1}\Big(\frac{m-1}{m}\sup_{x\in M^{n}}v\Big)^{\frac{\alpha -1}{m-1}}  \\
&\quad +\frac{(m-1)^{2}(k-1)}{m(m+1)}|\frac{2(m+1)}{(m-1)(k-1)}+1|\delta_{2},
\endaligned$$

$$\aligned
 H_{2}&=\frac{(m-1)(k-1)}{m(m+1)}|\frac{2(m+1)}{n-1}+(m-2\alpha+1)|\delta_{1}\Big(\frac{m-1}{m}\inf_{x\in M^{n}} v\Big)^{\frac{\alpha -1}{m-1}}  \\
& \quad+\frac{(m-1)^{2}(k-1)}{m(m+1)}|\frac{2(m+1)}{(m-1)(k-1)}+1|\delta_{2}.
\endaligned$$
Moreover, if $( M^{n}, g)$ has nonnegative Ricci curvature, letting $R\rightarrow \infty $, we have
following estimate for $\alpha \geq 1$,
\begin{equation}
\frac{|\nabla v|^{2}}{v}\leq C(m,n,\alpha,k,\delta_{1},\delta_{2},\sup_{x\in M^{n}}v),
\end{equation}
and for $\alpha < 1$
\begin{equation}
\frac{|\nabla v|^{2}}{v}\leq C'(m,n,\alpha,k,\delta_{1},\delta_{2},\sup_{x\in M^{n}}v,\inf_{x\in M^{n}} v).
\end{equation}
Where
$$C(m,n,\alpha,k,\delta_{1},\delta_{2},\sup_{x\in M^{n}}v)=\frac{C_4(m-1)}{m+1}\sup_{x\in M^n}v+H_1,$$
$$C'(m,n,\alpha,k,\delta_{1},\delta_{2},\sup_{x\in M^{n}}v,\inf_{x\in M^{n}}v)=\frac{C_4(m-1)}{m+1}\sup_{x\in M^n}v+H_2.$$

\end{thm}
\begin{corr}\label{coro1-1} Let $(M^{n}, g)$ be a noncompact complete Riemannian manifold without
 boundary. Suppose that $Ric_{V}^{k}(M^{n}) \geq 0$. Let $u(x)$ is a positive solution of the
equation \eqref{1Int1} , and $v=\frac{m}{m-1}u^{m-1}$.\\
If $\alpha \geq 1$, then
\begin{equation}
v(x)\leq v(y)\exp[r(x,y)\sqrt{\frac{C(m,n,\alpha,k,\delta_{1},\delta_{2},K,\sup_{x\in M^{n}}v)}{\inf_{x\in M^{n}}v}}].
\end{equation}
If $\alpha < 1$, then
\begin{equation}
v(x)\leq v(y)\exp[r(x,y)\sqrt{\frac{C'(m,n,\alpha,k,\delta_{1},\delta_{2},K,\sup_{x\in M^{n}}v,\inf_{x\in M^{n}}v)}{\inf_{x\in M^{n}}v}}].
\end{equation}

\end{corr}

\begin{thm}\label{Int-2} Let $(M^{n},g)$ be a complete Riemannian manifold without boundary. Suppose that $B_{2R}$ is a geodesic ball of radius $2R$ aroung $q\in M$ and $Ric_{V}^{k}(B_{2R}) \geq -K$
with $K \geq 0$. Let $u(x)$ is a positive solution of the equation \eqref{1Int1}. Let $v=\frac{m}{m-1}u^{m-1}$, $\beta>0$ , $|\nabla p|^{2} \leq |p|^{2}$ and $|\nabla \mu|^{2} \leq |\mu|^{2}$ . If $p\geq 0$, $\mu \geq 0$ and $\alpha\leq \frac{(k+1)(m+1)}{2(k-1)}$
, then we have
\begin{equation}
\sup_{B_q(R)}\frac{|\nabla v|^2}{v}\leq C_4[\frac{1}{R^{2}}(1+\sqrt{K}R)+2K+1]\sup_{x\in M^{n}}v,
\end{equation}
where $C_{4}$ is a constant depending on $k$.

Letting $R\rightarrow \infty $ , then we infer on a complete noncompact Riemannian manifold,
\begin{equation}
\frac{|\nabla v|^2}{v}\leq C_4[2K+1]\sup_{x\in M^{n}}v.
\end{equation}

\end{thm}
\begin{corr}\label{coro1-2}
Let $(M^{n},g)$ be a noncompact complete Riemannian manifold without boundary. Suppose that $Ric_{V}^{k}(M^{n}) \geq 0$ and $u(x)$ is a positive solution of the equation \eqref{1Int1}, where $p(x)$ and $\mu(x)$ are constant. If $p\geq 0$, $\mu\geq 0$ and $\alpha\leq\frac{(k+1)(m+1)}{2(k-1)}$, then $u$ is a constant.

\end{corr}

\begin{thm}\label{Int-3} Let $(M^{n},g)$ be a complete noncompact Riemannian manifold with
$Ric_{V}^{k}(M^{n}) \geq -K$ with $K \geq 0$. Let $u(x)$ is a positive solution to the equation
\begin{equation}
\Delta_{V}u^{m}+\mu u+p u^{\alpha}=0,
\end{equation}
where $p > 0$ and $\mu$ are constant. Let $v=\frac{m}{m-1}u^{m-1}$ and $1\leq \alpha \leq \frac{(m+1)(k+1)}{2(k-1)}$.

$( 1) $If $\mu\geq0$,
$$\aligned
p & \leq\frac{2(m-1)(k-1)K}{\frac{2(m+1)}{k-1}+(m-2\alpha+1)}(\frac{m}{m-1})^{\frac{\alpha-1}{m-1}}(\sup_{x\in M^{n}}v)^{\frac{m-\alpha}{m-1}} \\
&\quad-\frac{m-1}{\frac{2(m+1)}{k-1}+(m-2\alpha+1)}(\frac{2(m+1)}{(m-1)(k-1)}+1)\mu (\frac{m}{m-1})^{\frac{\alpha-1}{m-1}}(\frac{1}{\sup_{x\in M^{n}}v})^{\frac{\alpha-1}{m-1}}.
\endaligned$$
or $\mu<0$,
$$\aligned
p & \leq\frac{2(m-1)(k-1)K}{\frac{2(m+1)}{k-1}+(m-2\alpha+1)}(\frac{m}{m-1})^{\frac{\alpha-1}{m-1}}(\sup_{x\in M^{n}}v)^{\frac{m-\alpha}{m-1}} \\
&\quad-\frac{m-1}{\frac{2(m+1)}{k-1}+(m-2\alpha+1)}(\frac{2(m+1)}{(m-1)(k-1)}+1)\mu (\frac{m}{m-1})^{\frac{\alpha-1}{m-1}}(\frac{1}{\inf_{x\in M^{n}}v})^{\frac{\alpha-1}{m-1}},
\endaligned$$
then for any $x\in M^{n}$,
$$\aligned
\frac{|\nabla v|^{2}}{v}&\leq\frac{2(m-1)^{2}(k-1)^{2}}{m(m+1)}K\sup_{x\in M^{n}}v\\
&\quad-\frac{(m-1)(k-1)}{m(m+1)}[\frac{2(m+1)}{k-1}+(m-2\alpha+1)]p(\frac{m-1}{m}\sup_{x\in M^{n}}v)^{\frac{\alpha-1}{m-1}}\\
&\quad-\frac{(m-1)^{2}(k-1)}{m(m+1)}[\frac{2(m+1)}{(m-1)(k-1)}+1]\mu.
\endaligned$$

$( 2) $ If as $\mu\geq0$,
$$\aligned
p & \geq\frac{2(m-1)(k-1)K}{\frac{2(m+1)}{k-1}+(m-2\alpha+1)}(\frac{m}{m-1})^{\frac{\alpha-1}{m-1}}(\sup_{x\in M^{n}}v)^{\frac{m-\alpha}{m-1}} \\
&\quad-\frac{m-1}{\frac{2(m+1)}{k-1}+(m-2\alpha+1)}(\frac{2(m+1)}{(m-1)(k-1)}+1)\mu (\frac{m}{m-1})^{\frac{\alpha-1}{m-1}}(\frac{1}{\sup_{x\in M^{n}}v})^{\frac{\alpha-1}{m-1}}.
\endaligned$$
or $\mu<0$,
$$\aligned
p & \geq\frac{2(m-1)(k-1)K}{\frac{2(m+1)}{k-1}+(m-2\alpha+1)}(\frac{m}{m-1})^{\frac{\alpha-1}{m-1}}(\sup_{x\in M^{n}}v)^{\frac{m-\alpha}{m-1}} \\
&\quad-\frac{m-1}{\frac{2(m+1)}{k-1}+(m-2\alpha+1)}(\frac{2(m+1)}{(m-1)(k-1)}+1)\mu (\frac{m}{m-1})^{\frac{\alpha-1}{m-1}}(\frac{1}{\inf_{x\in M^{n}}v})^{\frac{\alpha-1}{m-1}},
\endaligned$$
then $v$ must be a constant.

\end{thm}

\section{Preliminaries}
\begin{lem}
Let $(M^{n},g)$ be a complete Rimannian mainfold without boundary. Assume that $B_{2R}$ is a geodesic ball of radius $2R$ around $q\in M$ and $Ric_{V}^{k}(B_{2R}) \geq -K$
with $K \geq 0$. Let $u(x)$ be a positive solution to the equation \eqref {1Int1} and $v=\frac{m}{m-1}u^{m-1}$. Let $\omega =\frac{|\nabla v|^2}{v} $ and $G = \varphi \omega$, where $\varphi(x)$ is a smooth cutoff function. Suppose that $G(x)$ reaches the maximum value at
$x_{0}$ and $\varphi(x_{0}) > 0$. Then at $x_{0}$,
\begin{equation}\begin{aligned}\label{2Int1}
\varphi\Delta_{V}\omega& \geq[\frac{k}{2(k-1)}+\frac{2m}{(k-1)(m-1)^{2}}+\frac{1}{m-1}]\frac{G^{2}}{v\varphi} \\
&\quad+[\frac{2m}{m-1}-\frac{k}{k-1}-\frac{2}{(k-1)(m-1)}]\frac{\nabla v\nabla \varphi}{v\varphi}G+\frac{k}{2(k-1)}\frac{|\nabla\varphi|^{2}}{\varphi^{2}}G \\
&\quad+\left[\big(\frac{2(m+1)}{(k-1)(m-1)}+1\big)\mu+\big(\frac{2(m+1)}{(k-1)(m-1)}+\frac{m-2\alpha+1}{m-1}\big) p(\frac{m-1}{m}v)^{\frac{\alpha-1}{m-1}}\right]\frac{G}{v} \\
&\quad-\frac{2}{v}\varphi|\nabla v|\Big(|\nabla \mu|+|\nabla p|(\frac{m-1}{m}v)^{\frac{\alpha-1}{m-1}}\Big)-\frac{2}{k-1}[\mu+ p(\frac{m-1}{m}v)^{\frac{\alpha-1}{m-1}}]\frac{\nabla v\nabla\varphi}{v} \\
&\quad+\frac2{k-1}[\mu+ p(\frac{m-1}{m}v)^{\frac{\alpha-1}{m-1}}]^{2}\frac{\varphi}{v}-2(k-1)KG.
\end{aligned}\end{equation}
\end{lem}
\proof Through a simple calculation, we can obtain
\begin{equation}\label{2Int2}
(m-1)v\Delta_{V} v+|\nabla v|^2=-p(m-1)\big(\frac{m-1}{m}\big)^{\frac{\alpha-1}{m-1}}v^{1+\frac{\alpha-1}{m-1}}-(m-1)\mu v.
\end{equation}
We can compute
\begin{equation}\begin{aligned}\label{2Int3}
\Delta_{V}\omega& =\Delta_{V}(\frac{|\nabla v|^{2}}{v}) \\
&=\frac{\Delta_{V}|\nabla v|^{2}}{v}-\frac{2\nabla|\nabla v|^{2}\nabla v}{v^{2}}-\frac{|\nabla v|^{2}\Delta_{V} v}{v^2}+\frac{2|\nabla v|^{4}}{v^{3}} \\
&=\frac{2}{v}[|\operatorname{Hess}v|^{2}+\nabla v\cdot\nabla\Delta_{V} v+\operatorname{Ric}_{V}(\nabla v,\nabla v)] \\
&\quad-\frac{2\nabla|\nabla v|^{2}\cdot\nabla v}{v^{2}}-\frac{|\nabla v|^{2}\Delta_{V} v}{v^{2}}+\frac{2|\nabla v|^{4}}{v^{3}}.
\end{aligned}\end{equation}
 Since $G$ attains its maximum at $x_{0}$, so $\nabla G = 0$. At $x_{0}$, we have
\begin{equation}\label{2Int4}
\nabla\omega=-\frac{G\nabla \varphi}{\varphi^{2}}
\end{equation}
and
\begin{equation}\label{2Int5}
\nabla|\nabla v|^{2}=-\frac{vG}{\varphi^{2}}\nabla \varphi+\frac{G}{\varphi}\nabla v.
\end{equation}
We choose an orthonormal frame $\{e_{1},\ldots,e_{n}\}$
around $x_{0}$, such that $|\nabla v|e_{1}=\nabla v$.
We have
\begin{equation}\label{2Int6}
\frac{|\nabla|\nabla v|^{2}|^{2}}{4|\nabla v|^{2}}=\sum_{j=1}^nv_{1j}^{2},
\end{equation}
\begin{equation}\label{2Int7}
\frac{\nabla v\cdot\nabla|\nabla v|^{2}}{2|\nabla v|^{2}}=v_{11}.
\end{equation}
By \eqref{2Int2},
$$\aligned
|\operatorname{Hess}v|^{2}&\geq v_{11}^{2}+2\sum_{\alpha=2}^{n}v_{1\alpha}^{2}+\sum_{\alpha=2}^{n}v_{\alpha\alpha}^{2} \\
&\geq v_{11}^{2}+2\sum_{\alpha=2}^nv_{1\alpha}^2+\frac{1}{n-1}(\sum_{\alpha=2}^nv_{\alpha\alpha})^{2} \\
&=v_{11}^{2}+2\sum_{\alpha=2}^nv_{1\alpha}^2+\frac{1}{n-1}(\Delta v - v_{11})^{2} \\
&=\frac{n}{n-1}v_{11}^2+2\sum_{\alpha=2}^nv_{1\alpha}^2-\frac{2v_{11}}{n-1}\Delta v+\frac{1}{n-1}(\Delta v)^{2} \\
&=\frac{n}{n-1}v_{11}^{2}+\frac{2v_{11}}{n-1}[\langle V,\nabla v\rangle+\frac{\omega}{m-1}+\mu+p(\frac{m-1}{m}v)^{\frac{\alpha-1}{m-1}}] \\
&\quad+\frac{1}{n-1}[\langle V,\nabla v\rangle+\frac{\omega}{m-1}+\mu+p(\frac{m-1}{m}v)^{\frac{\alpha-1}{m-1}}]^{2}+2\sum_{\alpha=2}^{n}v_{1\alpha}^{2}.
\endaligned$$
Then, we have
$$\aligned
|\operatorname{Hess}v|^{2}&+\mathrm{Ric}_{V}(\nabla v,\nabla v) \\
&\geq\frac{n}{n-1}v_{11}^{2}+\frac{2v_{11}}{n-1}[\langle V,\nabla v\rangle+\frac{\omega}{m-1}+\mu+p(\frac{m-1}{m}v)^{\frac{\alpha-1}{m-1}}] \\
&+\frac{1}{n-1}[\langle V,\nabla v\rangle+\frac{\omega}{m-1}+\mu+p(\frac{m-1}{m}v)^{\frac{\alpha-1}{m-1}}]^{2}+2\sum_{\alpha=2}^{n}v_{1\alpha}^{2} \\
&\quad-(k-1)K|\nabla v|^{2}+\frac{1}{k-n}\langle V,\nabla v\rangle^{2}\\
&=\frac{n}{n-1}v_{11}^2+2\sum_{\alpha=2}^nv_{1\alpha}^2+\frac{2v_{11}}{n-1}\big(\frac{\omega}{m-1}+\mu+p(\frac{m-1}{m}v)^{\frac{\alpha-1}{m-1}}\big)\\&\quad+\frac{1}{n-1}\big(\frac{\omega}{m-1}+\mu+p(\frac{m-1}{m}v)^{\frac{\alpha-1}{m-1}}\big)^2 +\frac{k-1}{(k-n)(n-1)}(\langle V,\nabla v\rangle)^{2}\\
&\quad+\frac{2\langle V,\nabla v\rangle}{n-1}\big(\frac{\omega}{m-1}+\mu+p(\frac{m-1}{m}v)^{\frac{\alpha-1}{m-1}}+v_{11}\big)-(k-1)K|\nabla v|^{2}.
\endaligned$$
Because of $a > 0$, we have $ax^{2}+bx\geq -\frac{b^{2}}{4a}$, so that
$$\aligned
&\frac{k-1}{(k-n)(n-1)}(\langle V,\nabla v\rangle)^{2}+\frac{2\langle V,\nabla v\rangle}{n-1}\big(\frac{\omega}{m-1}+\mu+p(\frac{m-1}{m}v)^{\frac{\alpha-1}{m-1}}+v_{11}\big)\\
&\geq-\frac{k-n}{(k-1)(n-1)}(\frac{\omega}{m-1}+\mu+p(\frac{m-1}{m}v)^{\frac{\alpha-1}{m-1}}+v_{11})^{2}.
\endaligned$$
Then,
$$\aligned
|\operatorname{Hess}v|^{2}&+\mathrm{Ric}_{V}(\nabla v,\nabla v) \\
&\geq
\frac{n}{n-1}v_{11}^2+2\sum_{\alpha=2}^nv_{1\alpha}^2+\frac{2v_{11}}{n-1}\big(\frac{\omega}{m-1}+\mu+p(\frac{m-1}{m}v)^{\frac{\alpha-1}{m-1}}\big)\\
&\quad+\frac{1}{n-1}\big(\frac{\omega}{m-1}+\mu+p(\frac{m-1}{m}v)^{\frac{\alpha-1}{m-1}}\big)^2-(k-1)K|\nabla v|^{2} \\
&\quad-\frac{k-n}{(k-1)(n-1)}(\frac{\omega}{m-1}+\mu+p(\frac{m-1}{m}v)^{\frac{\alpha-1}{m-1}}+v_{11})^{2}\\
&=\frac{k}{k-1}v_{11}^2+2\sum_{\alpha=2}^nv_{1\alpha}^2+\frac{2v_{11}}{k-1}(\frac{\omega}{m-1}+\mu+p(\frac{m-1}{m}v)^{\frac{\alpha-1}{m-1}})\\
&\quad+\frac{1}{k-1}\big(\frac{\omega}{m-1}+\mu+p(\frac{m-1}{m}v)^{\frac{\alpha-1}{m-1}}\big)^2-(k-1)K|\nabla v|^{2} \\
&\geq \frac{k}{k-1}\sum_{j=1}^nv_{1\j}^2+\frac{2v_{11}}{k-1}\big(\frac{\omega}{m-1}+\mu+p(\frac{m-1}{m}v)^{\frac{\alpha-1}{m-1}}\big)\\&\quad+\frac{1}{k-1}\big(\frac{\omega}{m-1}+\mu+p(\frac{m-1}{m}v)^{\frac{\alpha-1}{m-1}}\big)^2-(k-1)K|\nabla v|^{2}.
\endaligned$$
By\eqref{2Int4}, \eqref{2Int5}, \eqref{2Int6}and\eqref{2Int7}, we have
$$\aligned
|\operatorname{Hess}v|^{2}&+\mathrm{Ric}_{V}(\nabla v,\nabla v) \\
&\geq\frac{k}{k-1}\frac{|\nabla|\nabla v|^{2}|^{2}}{4|\nabla v|^{2}}+\frac{\nabla v\cdot\nabla|\nabla v|^{2}}{(k-1)|\nabla v|^{2}}[\frac{\omega}{m-1}+\mu+p(\frac{m-1}{m}v)^{\frac{\alpha-1}{m-1}}] \\
&\quad+\frac{1}{k-1}[\frac{\omega}{m-1}+\mu+p(\frac{m-1}{m}v)^{\frac{\alpha-1}{m-1}}]^{2}-(k-1)K|\nabla v|^{2} \\
&=\frac{k}{4(k-1)|\nabla v|^{2}}[\frac{G}{\varphi}\nabla v-\frac{vG}{\varphi^{2}}\nabla\varphi]^{2}+\frac{1}{k-1}[\frac{\omega}{m-1}+\mu+p(\frac{m-1}{m}v)^{\frac{\alpha-1}{m-1}}]^{2} \\
&\quad+\frac{1}{(k-1)|\nabla v|^{2}}\big[\frac{G}{\varphi}|\nabla v|^{2}-\frac{vG}{\varphi^{2}}\nabla v\nabla\varphi][\frac{\omega}{m-1}+\mu+p(\frac{m-1}{m}v)^{\frac{\alpha-1}{m-1}}] \\
&\quad-(k-1)K|\nabla v|^{2} \\
&=\frac{k}{4(k-1)}(\frac{G}{\varphi})^{2}+\frac{kv^{2}}{4(n-1)|\nabla v|^{2}}\frac{G^{2}}{\varphi^{4}}|\nabla\varphi|^{2}-\frac{k}{2(k-1)|\nabla v|^{2}}\cdot\frac{vG^{2}\nabla v\nabla\varphi}{\varphi^{3}} \\
&\quad+\frac{\omega^2}{(k-1)(m-1)^2}+\frac{2\omega}{(k-1)(m-1)}\big(\mu+p(\frac{m}{m-1}v)^{\frac{\alpha-1}{m-1}}\big) \\
&\quad+\frac{1}{k-1}\big(\mu+p\big(\frac{m}{m-1}v\big)^{\frac{\alpha-1}{m-1}}\big)^{2}+\frac{1}{(k-1)(m-1)}\frac{G}{\varphi}\omega+\frac{1}{k-1}\frac{G}{\varphi}(\mu+p\big(\frac{m}{m-1}v\big)^{\frac{\alpha-1}{m-1}}) \\
&\quad-\frac{1}{(k-1)|\nabla v|^{2}}\frac{vG}{\varphi^{2}}\nabla v\nabla\varphi\frac{w}{m-1}-\frac{1}{(k-1)|\nabla v|^{2}}\frac{vG}{\varphi^{2}}\nabla v\nabla\varphi\big(\mu+p\big(\frac{m}{m-1}v\big)^{\frac{\alpha-1}{m-1}}\big) \\
&\quad-(k-1)K|\nabla v|^{2}\\
&=\frac{k}{4(k-1)}(\frac{G}{\varphi})^{2}+\frac{kv}{4(k-1)}\frac{|\nabla\varphi|^{2}}{\varphi^{2}}\frac{G}{\varphi}-\frac{k}{2(k-1)}\frac{G}{\varphi}\frac{\nabla v\nabla\varphi}{\varphi} \\
&\quad+\frac{1}{(k-1)(m-1)^{2}}(\frac{G}{\varphi})^{2}+\frac{1}{k-1}\big(\mu+p(\frac{m-1}{m}v)^{\frac{\alpha-1}{m-1}}\big)^{2}
\\&\quad+\frac{2}{(k-1)(m-1)}\big(\mu+p(\frac{m-1}{m}v)^{\frac{\alpha-1}{m-1}}\big)\frac{G}{\varphi} +\frac{1}{(k-1)(m-1)}(\frac{G}{\varphi})^2 \\
&\quad+\frac1{k-1}\big(\mu+p(\frac{m-1}{m}v)^{\frac{\alpha-1}{m-1}}\big)\frac G\varphi-\frac1{(k-1)(m-1)}\frac G\varphi\frac{\nabla v\nabla\varphi}\varphi \\
&\quad-\frac{1}{k-1}\big(\mu+p(\frac{m-1}{m}v)^{\frac{\alpha-1}{m-1}}\big)\frac{\nabla v\nabla\varphi}{\varphi}-(k-1)K|\nabla v|^{2} \\
&=[\frac{k}{4(k-1)}+\frac{m}{(k-1)(m-1)^{2}}]\frac{G^{2}}{\varphi^{2}}-[\frac{k}{2(k-1)}+\frac{1}{(k-1)(m-1)}]\frac{G}{\varphi^{2}}\nabla v\nabla\varphi \\
&\quad+\frac{k}{4(k-1)}\frac{|\nabla\varphi|^{2}}{\varphi^{3}}vG+\frac{m+1}{(k-1)(m-1)}\big(\mu+p(\frac{m-1}{m}v)^{\frac{\alpha-1}{m-1}}\big)\frac{G}{\varphi} \\
&\quad-\frac{1}{k-1}\big(\mu+p(\frac{m-1}{m}v)^{\frac{\alpha-1}{m-1}}\big)\frac{\nabla v\nabla\varphi}{\varphi}+\frac{1}{k-1}\big(\mu+p(\frac{m-1}{m}v)^{\frac{\alpha-1}{m-1}}\big)^{2} \\
&\quad-(k-1)K|\nabla v|^{2}.
\endaligned$$
Noting that $v > 0$ when $m > 1$.Substituting the inequality into \eqref{2Int3}, then we can get \eqref{2Int1}.
\endproof

\section{The proof of main theorems.}
In this section, we prove Theorem \ref{Int-1}.
\proof
Define the smooth cutoff function $\varphi : M \rightarrow \mathbb{R}$ by $\varphi(x) = \theta(\frac{r(x)}{R})$. Suppose that
$G = \varphi \omega = \varphi \frac{|\nabla v|^{2}}{v}$  and attains its maximal value at $x_{0}\in B_{2R}$. We assume $G(x_{0}) > 0$, then $\varphi(x_{0})> 0$. Then at $x_{0}$, we have
$$\aligned
\Delta_{V} G&=\Delta_{V}\varphi\cdot \omega+2\nabla\varphi\nabla \omega+\varphi\Delta_{V} \omega\\
&=\Delta_{V}\varphi\cdot \omega-2G\frac{|\nabla\varphi|^{2}}{\varphi^{2}}+\varphi\Delta_{V} \omega\\
&=\frac{\Delta_{V}\varphi}{\varphi}G-2G\frac{|\nabla\varphi|^2}{\varphi^2}+\varphi\Delta_{V} \omega.
\endaligned$$
Note that,
$$\nabla\varphi=\frac{\theta'\nabla r}{R},$$
by the weighted Laplacian comparison theorem,
$$\Delta_{V}\varphi=\frac{\theta''}{R^2}+\frac{\theta'\Delta_{V} r}{R}\geq\frac{\theta''}{R^2}+\frac{(k-1)(1+\sqrt{K}R)\theta'}{R^2}.$$
Since $\Delta_{V} G \leq 0$ at $x_{0}$, we have
\begin{equation}\begin{aligned}\label{3Int1}
\text{0}& \geq [\frac{\theta^{\prime\prime}}{\theta R^{2}}+\frac{(k-1)(1+\sqrt{K}R)\theta^{\prime}}{\theta R^{2}}]G-\frac{3k-4}{2(k-1)}\frac{(\theta^{\prime})^{2}}{R^{2}\theta^{2}}G \\
&\quad+\frac{(m-1)(km+k-2)+4m}{2(m-1)^2(k-1)}\frac{G^{2}}{v\theta}-\frac{(m+1)(k-2)}{(k-1)(m-1)}\frac{G\sqrt{G}|\theta^{\prime}|}{R\theta\sqrt{v\theta}} \\
&\quad+\left[(\frac{2(m+1)}{(k-1)(m-1)}+1)\mu+(\frac{2(m+1)}{(m-1)(k-1)}+\frac{m-2\alpha+1}{m-1})p(\frac{m-1}{m}v)^{\frac{\alpha-1}{m-1}}\right]\frac{G}{v} \\
&\quad-2(k-1)KG-\frac{2}{v}\theta|\nabla v|\big(|\nabla \mu|+|\nabla p|(\frac{m}{m-1}v)^{\frac{\alpha-1}{m-1}}\big) \\
&\quad-\frac{2}{k-1}\big(\mu+p(\frac{m-1}{m}v)^{\frac{\alpha-1}{m-1}}\big)\frac{\sqrt{G}|\theta^{\prime}|}{R\sqrt{v\theta}}+\frac{2}{k-1}\big(\mu+p(\frac{m-1}{m}v)^{\frac{\alpha-1}{m-1}}\big)^{2}\frac{\theta}{v}.
\end{aligned}\end{equation}
Bacause $\mu p \geq 0$, $|\nabla p|^{2} \leq |p|^{2}$ and $|\nabla \mu|^{2} \leq |\mu|^{2}$, we have
$$\frac{1}{k-1}\big(\mu+p(\frac{m-1}{m}v)^{\frac{\alpha-1}{m-1}}\big)^{2}\frac{\theta}{v}=\frac{1}{k-1}\big(|\mu|+|p|(\frac{m-1}{m}v)^{\frac{\alpha-1}{m-1}}\big)^{2}\frac{\theta}{v}$$
$$\aligned-\frac{2}{v}\theta|\nabla v|\big(|\nabla \mu|+|\nabla p|(\frac{m}{m-1}v)^{\frac{\alpha-1}{m-1}}\geq -\frac{2}{v}\theta|\nabla v|\big(|\mu|+| p|(\frac{m}{m-1}v)^{\frac{\alpha-1}{m-1}}\endaligned$$
Using the inequality $ax^{2}+bx\geq -\frac{b^{2}}{4a}$, $a>0$, we have
\begin{equation}\begin{aligned}\label{3Int2}
&\frac{1}{k-1}\big(\mu+p(\frac{m-1}{m}v)^{\frac{\alpha-1}{m-1}}\big)^{2}\frac{\theta}{v}-\frac{2}{k-1}(\mu+p(\frac{m-1}{m}v)^{\frac{\alpha-1}{m-1}})\frac{\sqrt{G}|\theta^{\prime}|}{R\sqrt{v\theta}}\\
&\geq -\frac{G(\theta')^{2}}{(k-1)R^{2}\theta^{2}},
\end{aligned}\end{equation}
and
\begin{equation}\begin{aligned}\label{3Int3}
&\frac{1}{k-1}\big(|\mu|+|p|(\frac{m-1}{m}v)^{\frac{\alpha-1}{m-1}}\big)^{2}\frac{\theta}{v}-\frac{2}{v}\theta|\nabla v|\big(| \mu|+| p|(\frac{m}{m-1}v)^{\frac{\alpha-1}{m-1}}\big) \\&\geq-(k-1)G.
\end{aligned}\end{equation}
By the Cauchy inequality, it follows that
\begin{equation}\begin{aligned}\label{3Int4}
-\frac{G\sqrt{G}|\theta'|}{R\theta\sqrt{v\theta}}\geq-\frac{G^2}{2v\theta}-\frac{G(\theta')^2}{2R^2\theta^2}.
\end{aligned}\end{equation}
Substituting \eqref{3Int2}, \eqref{3Int3}, and\eqref{3Int4}to\eqref{3Int1}, we have
\begin{equation}\begin{aligned}\label{3Int5}
0& \geq[\frac{\theta^{\prime\prime}}{\theta R^{2}}+\frac{(k-1)(1+\sqrt{K}R)\theta^{\prime}}{\theta R^{2}}]G-\frac{3k-2}{2(k-1)}\frac{(\theta^{\prime})^{2}}{R^{2}\theta^{2}}G \\
&\quad+\frac{m(m+1)}{(m-1)^2(k-1)}\frac{G^2}{v\theta}-\frac{(m+1)(k-2)}{2(m-1)(k-1)}\frac{G(\theta')^2}{R^2\theta^2} \\
&\quad-2(k-1)KG-(k-1)G \\
&\quad+[(\frac{2(m+1)}{(k-1)(m-1)}+1)\mu+(\frac{2(m+1)}{(k-1)(m-1)}+\frac{m-2\alpha+1}{m-1})p(\frac{m-1}{m}v)^{\frac{\alpha-1}{m-1}}]\frac{G}{v}.
\end{aligned}\end{equation}
Next, we construct a smooth function $\theta(t) : [0, +\infty) \rightarrow [0, 1]$
$$\theta(t)=\begin{cases}1,&0\le t\le1\\0,&t>2\end{cases},$$
and
\begin{equation}\label{3Int6}
-C_1\sqrt{\theta}\leq\theta'\leq 0,\quad|\theta''|\leq C_2\theta.
\end{equation}
Use \eqref{3Int6},
\begin{equation}\begin{aligned}\label{3Int9}
0&\geq-[\frac{C_{2}}{R^{2}}+\frac{(k-1)(1+\sqrt{K}R)C_{1}}{\sqrt{\theta}R^{2}}]G+\frac{m(m+1)}{(m-1)^{2}(k-1)}\frac{G^{2}}{v\theta} \\
&\quad-[\frac{2m}{m-1}-\frac{k}{(m-1)(k-1)}]\frac{C_{1}^{2}}{R^{2}\theta}G-2(k-1)KG-(k-1)G \\
&\quad+[(\frac{2(m+1)}{(k-1)(m-1)}+1)\mu+(\frac{2(m+1)}{(k-1)(n-1)}+\frac{m-2\alpha+1}{m-1})p(\frac{m-1}{m}v)^{\frac{\alpha-1}{m-1}}]\frac{G}{v} \\
&\geq-[\frac{C_{2}}{R^{2}}+\frac{(k-1)(1+\sqrt{K}R)C_{1}}{\sqrt{\theta}R^{2}}]G+\frac{m(m+1)}{(m-1)^{2}(k-1)}\frac{G^{2}}{v\theta} \\
&\quad-\frac{2m}{m-1}\frac{C_1^2}{R^2\theta}G-2(k-1)KG-(k-1)G \\
&\quad+[(\frac{2(m+1)}{(k-1)(m-1)}+1)\mu+(\frac{2(m+1)}{(k-1)(m-1)}+\frac{m-2\alpha+1}{m-1})p(\frac{m-1}{m}v)^{\frac{\alpha-1}{m-1}}]\frac{G}{v}.
\end{aligned}\end{equation}
Multiply both sides of \eqref{3Int9} by $v\theta $ , and  $0 \leq \theta \leq 1$ we have for $\alpha \geq 1$
\begin{equation}\begin{aligned}\label{3Int10}
0& \geq\frac{m(m+1)}{(m-1)^{2}(k-1)}G^{2}-\frac{2m}{m-1}\frac{C_{1}^{2}}{R^{2}}G\sup_{x\in M^{n}}v \\
&\quad-[\frac{C_2}{R^2}+\frac{(k-1)(1+\sqrt{K}R)C_1}{R^2}]G\sup_{x\in M^{n}}v \\
&\quad-[2(k-1)K+(k-1)]G\sup_{x\in M^{n}}v \\
&\quad-|\frac{2(m+1)}{(k-1)(m-1)}+1|\cdot|\mu|G \\ &\quad-|\frac{2(m+1)}{(k-1)(m-1)}+\frac{m-2\alpha+1}{m-1}|\cdot|p|(\frac{m-1}{m}\sup_{x\in M^{n}}v)^{\frac{\alpha-1}{m-1}}G.
\end{aligned}\end{equation}
Also, for $\alpha < 1$ we have
\begin{equation}\begin{aligned}\label{3Int11}
0& \geq\frac{m(m+1)}{(m-1)^{2}(k-1)}G^{2}-\frac{2m}{m-1}\frac{C_{1}^{2}}{R^{2}}G\sup_{x\in M^{n}}v \\
&\quad-[\frac{C_2}{R^2}+\frac{(k-1)(1+\sqrt{K}R)C_1}{R^2}]G\sup_{x\in M^{n}}v \\
&\quad-[2(k-1)K+(k-1)]G\sup_{x\in M^{n}}v \\
&\quad-|\frac{2(m+1)}{(k-1)(m-1)}+1|\cdot|\mu|G\\
&\quad-|\frac{2(m+1)}{(k-1)(m-1)}+\frac{m-2\alpha+1}{m-1}|\cdot|p|(\frac{m-1}{m}\inf_{x\in M^{n}}v)^{\frac{\alpha-1}{m-1}}G.
\end{aligned}\end{equation}
We observe that
\begin{equation}\label{3Int12}
\left[\frac{C_2}{R^2}+\frac{(k-1)(1+\sqrt{K}R)C_1}{R^2}\right]\sup\limits_{x\in M^n}v\leq\frac{C_3}{R^2}(1+\sqrt{K}R)\sup\limits_{x\in M^n}v,
\end{equation}
for some constant $C_{3}$ depending only on $k$.

On the contrary, for the inequality  $Ax^{2}-Bx\leq 0$ with $A > 0$, $B > 0$, we have $x\leq \frac{B}{A}$.
 By applying this result to equations  \eqref{3Int10} and \eqref{3Int11}, and  considering condition \eqref{3Int12}, we find the maximum point $x_{0}$ for $\alpha \geq 1$
$$\aligned
&\sup_{B_{q}(R)}\omega(x) \leq \varphi \omega(x_{0})=G(x_{0})\\&\leq\frac{2(m-1)(k-1)}{m+1}\frac{C_{1}^{2}}{R^{2}}\sup_{x\in M^{n}}v+\frac{(m-1)^{2}(k-1)}{m(m+1)}\frac{C_{3}}{R^{2}}(1+\sqrt{K}R)\sup_{x\in M^{n}}v \\
&\quad+\frac{(m-1)^{2}(k-1)^{2}}{m(m+1)}(2K+1)\sup_{x\in M^{n}}v \\
&\quad+\frac{(m-1)^{2}(k-1)}{m(m+1)}|\frac{2(m+1)}{(k-1)(m-1)}+1|\cdot|\mu|\\
&\quad+\frac{(m-1)(k-1)}{m(m+1)}|\frac{2(m+1)}{k-1}+(m-2\alpha+1)|\cdot|p|(\frac{m-1}{m}\sup_{x\in M^{n}}v)^{\frac{\alpha-1}{m-1}},
\endaligned$$
anf for $\alpha < 1$
$$\aligned
&\sup_{B_{q}(R)}\omega(x) \leq \varphi \omega(x_{0})=G(x_{0})\\&\leq\frac{2(m-1)(k-1)}{m+1}\frac{C_{1}^{2}}{R^{2}}\sup_{x\in M^{n}}v+\frac{(m-1)^{2}(k-1)}{m(m+1)}\frac{C_{3}}{R^{2}}(1+\sqrt{K}R)\sup_{x\in M^{n}}v \\
&\quad+\frac{(m-1)^{2}(k-1)^{2}}{m(m+1)}(2K+1)\sup_{x\in M^{n}}v \\
&\quad+\frac{(m-1)^{2}(k-1)}{m(m+1)}|\frac{2(m+1)}{(k-1)(m-1)}+1|\cdot|\mu|\\
&\quad+\frac{(m-1)(k-1)}{m(m+1)}|\frac{2(m+1)}{k-1}+(m-2\alpha+1)|\cdot|p|(\frac{m-1}{m}\inf_{x\in M^{n}}v)^{\frac{\alpha-1}{m-1}},
\endaligned$$
The proof is concluded.
\endproof

We prove Theorem \ref{Int-2}
\proof
A simple calculation shows that $\frac{2(m+1)}{(k-1)(m-1)}+\frac{m-2\alpha+1}{m-1}\geq 0 $
as $p\geq 0$ and $\mu\geq 0$. Hence, discarding the last nonnegative term in \eqref{3Int9}, we have the required inequality
$$\aligned
0\geq & -[\frac{C_{2}}{R^{2}}+\frac{(k-1)(1+\sqrt{KR})C_{1}}{\sqrt{\theta}R^{2}}]G+\frac{m(m+1)}{(m-1)^{2}(k-1)}\frac{G^{2}}{v\theta} \\
&- \frac{2m}{m-1}\frac{C_{1}^{2}}{R^{2}\theta}G-(k-1)(2K+1)G.
\endaligned$$
Multiply both sides of \eqref{3Int9} by $v\theta $, and  $0 \leq \theta \leq 1$, we have
$$\aligned
0&\geq\frac{m(m+1)}{(m-1)^{2}(k-1)}G^{2}-\frac{2m}{m-1}\frac{C_{1}^{2}}{R^{2}}G\sup_{x\in M^{n}}v \\
&\quad-[\frac{C_2}{R^{2}}+\frac{(k-1)(1+\sqrt{K}R)C_1}{R^{2}}]G \sup_{x\in M^{n}} v\\
&\quad-(k-1)(2K+1)G \sup_{x\in M^{n}}v.
\endaligned$$
Hence, at the the maximum point $x_{0}$ we obtain
$$\aligned
\sup_{B_{q}(R)}\omega(x)& \leq\varphi \omega(x_{0})=G(x_{0}) \\
&\leq\frac{2(m-1)(k-1)}{m+1}\frac{C_{1}^{2}}{R^{2}}\sup_{x\in M^{n}}v\\
&\quad+\frac{(m-1)^{2}(k-1)}{m(m+1)}\frac{C_{3}}{R^{2}}(1+\sqrt{K}R)\sup_{x\in M^{n}}v \\
&\quad+\frac{(m-1)^{2}(k-1)^{2}}{m(m+1)}(2K+1)\sup_{x\in M^{n}}v.
\endaligned$$
We used \eqref{3Int12} in the above inequality. The proof is concluded.
\endproof Next, we prove Theorem \ref{Int-3}
\proof It is easy to see $\frac{2(m+1)}{(k-1)(m-1)}+\frac{m-2\alpha+1}{m-1}\geq 0 $ for
$1\leq \alpha\leq\frac{(k+1)(m+1)}{2(k-1)}$ . Then we have form \eqref{3Int9},

\begin{equation}\begin{aligned}\label{3Int13}
0& \geq\frac{m(m+1)}{(m-1)^{2}(k-1)}G^{2}-\frac{2m}{m-1}\frac{C_{1}^{2}}{R^{2}} \sup_{x\in M^{n}}vG \\
&\quad-\Big[\frac{C_2}{R^2}+\frac{(k-1)(1+\sqrt{K}R)C_1}{R^2}]G\sup_{x\in M^{n}}v \\
&\quad-[2(k-1)K+(k-1)]G\sup_{x\in M^{n}}v \\
&\quad+[\frac{2(m+1)}{(k-1)(m-1)}+1]\mu G\\
&\quad+[\frac{2(m+1)}{(k-1)(m-1)}+\frac{m-2\alpha+1}{m-1}]p(\frac{m-1}{m}\sup_{x\in M^{n}}v)^{\frac{\alpha-1}{m-1}}G.
\end{aligned}\end{equation}
By \eqref{3Int13}, and \eqref{3Int12} we attain at the the maximum point $x_{0}$,
$$\aligned
\sup_{B_{q}(R)}\omega(x)& \leq\varphi \omega(x_{0})=G(x_{0}) \\
&\leq\frac{2(m-1)(k-1)}{m+1}\frac{C_{1}^{2}}{R^{2}}\sup_{x\in M^{n}}v+\frac{(m-1)^{2}(k-1)}{m(m+1)}\frac{C_{3}}{R^{2}}(1+\sqrt{K}R)\sup_{x\in M^{n}}v \\
&\quad+\frac{2(m-1)^{2}(k-1)^{2}}{m(m+1)}K\sup_{x\in M^{n}}v \\
&\quad-\frac{(m-1)^{2}(k-1)}{m(m+1)}[\frac{2(m+1)}{(k-1)(m-1)}+1]\mu\\
&\quad-\frac{(m-1)(k-1)}{m(m+1)}[\frac{2(m+1)}{k-1}+(m-2\alpha+1)]p(\frac{m-1}{m}\sup_{x\in M^{n}}v)^{\frac{\alpha-1}{m-1}},
\endaligned$$
Letting $R\rightarrow \infty $ , when $\mu\geq0$
$$\aligned
p & \leq\frac{2(m-1)(k-1)K}{\frac{2(m+1)}{k-1}+(m-2\alpha+1)}(\frac{m}{m-1})^{\frac{\alpha-1}{m-1}}(\sup_{x\in M^{n}}v)^{\frac{m-\alpha}{m-1}} \\
&\quad-\frac{m-1}{\frac{2(m+1)}{k-1}+(m-2\alpha+1)}(\frac{2(m+1)}{(m-1)(k-1)}+1)\mu (\frac{m}{m-1})^{\frac{\alpha-1}{m-1}}(\frac{1}{\sup_{x\in M^{n}}v})^{\frac{\alpha-1}{m-1}},
\endaligned$$
when $\mu<0$
$$\aligned
p & \leq\frac{2(m-1)(k-1)K}{\frac{2(m+1)}{k-1}+(m-2\alpha+1)}(\frac{m}{m-1})^{\frac{\alpha-1}{m-1}}(\sup_{x\in M^{n}}v)^{\frac{m-\alpha}{m-1}} \\
&\quad-\frac{m-1}{\frac{2(m+1)}{k-1}+(m-2\alpha+1)}(\frac{2(m+1)}{(m-1)(k-1)}+1)\mu (\frac{m}{m-1})^{\frac{\alpha-1}{m-1}}(\frac{1}{\inf_{x\in M^{n}}v})^{\frac{\alpha-1}{m-1}},
\endaligned$$
and

$$\aligned
\frac{|\nabla v|^{2}}{v}&\leq\frac{2(m-1)^{2}(k-1)^{2}}{m(m+1)}K\sup_{x\in M^{n}}v\\
&\quad-\frac{(m-1)(k-1)}{m(m+1)}[\frac{2(m+1)}{k-1}+(m-2\alpha+1)]p(\frac{m-1}{m}\sup_{x\in M^{n}}v)^{\frac{\alpha-1}{m-1}}\\
&\quad-\frac{(m-1)^{2}(k-1)}{m(m+1)}[\frac{2(m+1)}{(m-1)(k-1)}+1]\mu.
\endaligned$$
On the contrary, as
$\mu\geq0$
$$\aligned
p & \geq\frac{2(m-1)(k-1)K}{\frac{2(m+1)}{k-1}+(m-2\alpha+1)}(\frac{m}{m-1})^{\frac{\alpha-1}{m-1}}(\sup_{x\in M^{n}}v)^{\frac{m-\alpha}{m-1}} \\
&\quad-\frac{m-1}{\frac{2(m+1)}{k-1}+(m-2\alpha+1)}(\frac{2(m+1)}{(m-1)(k-1)}+1)\mu (\frac{m}{m-1})^{\frac{\alpha-1}{m-1}}(\frac{1}{\sup_{x\in M^{n}}v})^{\frac{\alpha-1}{m-1}},
\endaligned$$
when $\mu<0$
$$\aligned
p & \geq\frac{2(m-1)(k-1)K}{\frac{2(m+1)}{k-1}+(m-2\alpha+1)}(\frac{m}{m-1})^{\frac{\alpha-1}{m-1}}(\sup_{x\in M^{n}}v)^{\frac{m-\alpha}{m-1}} \\
&\quad-\frac{m-1}{\frac{2(m+1)}{k-1}+(m-2\alpha+1)}(\frac{2(m+1)}{(m-1)(k-1)}+1)\mu (\frac{m}{m-1})^{\frac{\alpha-1}{m-1}}(\frac{1}{\inf_{x\in M^{n}}v})^{\frac{\alpha-1}{m-1}},
\endaligned$$
we derive that $v$ must be constant.
\endproof

Proof of Corollary \ref{coro1-1}. \proof Let minimal geodesic $\gamma(s) : [0, 1] \rightarrow M^{n}$, so that $\gamma(0) = y$, $\gamma(1) = x$, then
$$\aligned
\ln{\frac{v(x)}{v(y)}}& =\int_{0}^{1}\frac{d\ln(v(\gamma(s)))}{ds}=\int_{0}^{1}\frac{\nabla v\cdot\gamma^{\prime}}{v(\gamma(s))}ds \\
&\leq\int_0^1\frac{|\nabla v|\cdot|\gamma^{\prime}|}{|v(\gamma(s))|}ds=r(x,y)\int_0^1\frac{|\nabla v|}{|v(\gamma(s))|}ds \\
&\leq r(x,y)\int_{0}^{1}\sqrt{\frac{C(m,n,\alpha,k,\delta_{1},\delta_{2},K,\sup_{x\in M^{n}}v)}{\inf_{x\in M^{n}}v}}ds\\
&=r(x,y)\sqrt{\frac{C(m,n,\alpha,k,\delta_{1},\delta_{2},K,\sup_{x\in M^{n}}v)}{\inf_{x\in M^n}v}}.
\endaligned$$
\endproof

\textbf{Acknowledgments}.
We are grateful to Professor Yingbo Han for his encouragement and concern. We also thank
Professor Fanqi Zeng for his valuable comments and answers to questions on the paper.

\bibliographystyle{Plain}

\end{document}